%
%
\magnification=\magstep1
\hsize=6.5truein
\vsize=9.0truein
\sfcode`"=0
\input epsf
\font\sc=cmcsc10

\font\bbf=cmbx10 scaled\magstep1
\font\rms=cmr10 scaled840
\input epsf

\def\DASIM{1.1}
\def\LORINTRO{1.2}
\def\LORINTROS{1.3}
\def\XDISCS{Figure$\ $1.1}
\def\TWODNS{1.4}
\def\PLAMBDAD{1.5}
\def\LORENZ{2.1}

\def\POINL{2.2}
\def\SYMML{2.3}
\def\NONLINL{2.4}
\def\BOUNDK{Theorem$\ $2.2}
\def\ETWO{2.5}
\def\YOUNGS{Theorem$\ $2.3}
\def\DIFEQL{2.6}
\def\EXGROWTHL{Lemma$\ $2.4}
\def\EXGA{2.7}
\def\EONE{2.8}
\def\LORENZDA{Theorem$\ $2.5}
\def\LORDE{2.9}
\def\THISB{2.10}
\def\THISA{2.11}
\def\THISC{2.12}

\def\UBOUNDB{2.14}
\def\ANBND{Corollary$\ $2.6}
\def\VARHL{Corollary$\ $2.7}

\def\BNDSLAVE{Theorem$\ $2.8}

\def\INDUL{2.17}
\def\UBOUND{Corollary$\ $2.9}
\def\FTWODNS{3.1}
\def\DEFAB{3.2}
\def\STRONGSOL{Theorem$\ $3.1}
\def\POINCARE{Theorem$\ $3.2}
\def\PYTHAG{Theorem$\ $3.3}
\def\ORTHAA{3.3}
\def\ORTHAB{3.4}
\def\SOBOLEV{Theorem$\ $3.4}
\def\NONLIN{Theorem$\ $3.5}
\def\APRIORI{Theorem$\ $3.6}
\def\EFOUR{3.5}
\def\DIFEQ{3.6}
\def\NOINDUCT{Theorem$\ $3.7}

\def\EXGROWTH{Theorem$\ $3.8}
\def\NSDA{Theorem$\ $3.9}

\def\DEFGEE{3.9}
\def\CONETWO{3.10}
\def\DEFEMM{3.11}
\def\BIGH{Corollary$\ $3.10}
\def\ZEROETA{Corollary$\ $3.11}

\def\UNTA{Theorem$\ $3.13}

\def\CONFOI{2}
\def\CONFOITEM{3}
\def\FJKT{4}
\def\FOIPRO{5}
\def\KEVIN{6}
\def\JONTITA{7}
\def\JONTITB{8}
\def\TAYLOR{9}
\def\LADY{10}
\def\OLSTIT{11}
\def\OLSTITB{12}

\def\PECAR{14}
\def\ROB{15}
\def\TEMA{16}
\def\TEMB{17}
\def\WINGARD{18}
\def\WU{19}
\def\YANG{20}

\newread\jobnameaux
\openin\jobnameaux=\jobname.aux
\ifeof\jobnameaux
    \closein\jobnameaux
    \message{...}
    \message{Reference tags being generated please ignore any undefined}
    \message{control sequences and run TeX one more time}
    \message{...}
\else
    \closein\jobnameaux
    \relax\input\jobname.aux
\fi
{\catcode`\@=0\catcode`\\=12@gdef@toauxs#1#2{%
@immediate@write@xrefs{\def\#1{#2}}}}
\def\toaux#1#2{\toauxs{#1}{#2}} 
\newwrite\xrefs\immediate\openout\xrefs=\jobname.aux
\newcount\chano\newcount\eqano\newcount\refno
\newcount\itmno\newcount\secno
\def\ref#1{\global\advance\refno by1\indent
\llap{[\the\refno]\enspace}\hang\toaux{#1}{\the\refno}\ignorespaces}
\def\qed{\vbox{\hrule\hbox{\vrule height6pt\hskip6pt\vrule}\hrule}}
\def\ifwrt#1#2{\if!#1!\else{\toaux{#1}{#2}}\fi}
\def\ses#1#2{\medskip\noindent{\bf #2.}\ifwrt{#1}{#2}}
\def\sesn#1#2{\noindent{\bf #2.}\ifwrt{#1}{#2}}
\def\dfn#1{\global\advance\itmno
    by1\ses{#1}{Definition$\ $\the\chano.\the\itmno}}
\def\dfnn#1{\global\advance\itmno
    by1\sesn{#1}{Definition$\ $\the\chano.\the\itmno}}

\def\fig#1{\global\advance\itmno
    by1\sesn{#1}{Figure$\ $\the\chano.\the\itmno}}
\def\thm#1{\global\advance\itmno
    by1\ses{#1}{Theorem$\ $\the\chano.\the\itmno}}
\def\thmn#1{\global\advance\itmno
    by1\sesn{#1}{Theorem$\ $\the\chano.\the\itmno}}
\def\prp#1{\global\advance\itmno
    by1\ses{#1}{Proposition$\ $\the\chano.\the\itmno}}
\def\rmk#1{\global\advance\itmno
    by1\ses{#1}{Remark$\ $\the\chano.\the\itmno}}
\def\lem#1{\global\advance\itmno
    by1\ses{#1}{Lemma$\ $\the\chano.\the\itmno}}
\def\lemn#1{\global\advance\itmno
    by1\sesn{#1}{Lemma$\ $\the\chano.\the\itmno}}
\def\cor#1{\global\advance\itmno
    by1\ses{#1}{Corollary$\ $\the\chano.\the\itmno}}
\def\corn#1{\global\advance\itmno
    by1\sesn{#1}{Corollary$\ $\the\chano.\the\itmno}}
\def\fct#1{\global\advance\itmno
    by1\ses{#1}{Fact$\ $\the\chano.\the\itmno}}
\def\beginsection#1{\global\eqano=0\global
    \itmno=0\global\advance\chano by1\bigskip
    \noindent{\bf \the\chano.  #1}\par\nobreak\medskip
    \noindent\ignorespaces}
\def\beginsubsection#1{
    \global\advance\secno by1\bigskip
    \noindent{\bf \the\chano.\the\secno. #1}\par\nobreak\medskip
    \noindent\ignorespaces}
\def\eqn#1{\global\advance\eqano by1\ifwrt{#1}{\the\chano.\the\eqano}%
    \hbox{\the\chano.\the\eqano}}
\def\edfn{\medskip\noindent}

\def\ethm{\medskip\noindent}

\def\ecor{\medskip\noindent}

\def\prf{\smallskip\noindent{\bf Proof: }}

\def\eprfn{\hfill{\qed}\par}
\def\eprf{\eprfn\medskip\noindent}

\headline{{\rm \hfill 
Tue Oct 12 04:02:13 MDT 2010
Version K22}}

\def\R{{\bf R}}

\def\N{{\bf N}}

\def\V{{\cal V}}

\def\J{{\cal J}}

\tolerance=3000

\vglue0pt
\bigskip
\vfill
{
\baselineskip=23pt
\centerline{\bbf Discrete Data Assimilation in the Lorenz}
\centerline{\bbf and 2D Navier--Stokes Equations}
}
\bigskip
\centerline{
{\sc Kevin Hayden}\footnote{$^1$}{{\rms Department of
Mathematics, Northern Arizona University, Flagstaff AZ 86011, USA}},
{\sc Eric Olson}\footnote{$^2$}{{\rms Department of Mathematics,
University of Nevada, Reno NV 89557, USA}}
and {\sc Edriss S.~Titi}\footnote{$^3$}{{\rms
Department of
Mathematics, University of California, Irvine CA
92697, USA}}$^,$\footnote{$^4$}{{\rms Mechanical
and Aerospace Engineering, University
of California, Irvine CA 92697, USA}}%
$^{,}$\footnote{$^5$}{{\rms
Department of Computer Science and Applied Mathematics,
Weizmann Institute of Science, Rehovot 76100, Israel}}
}
\bigskip

\def\Z{{\bf Z}}
\def\J{{\cal J}}

\vglue0pt
\medskip
{
\narrower
\narrower
\noindent
Consider a continuous dynamical system for which partial
information about its current state is observed at a sequence 
of discrete times.
Discrete data assimilation inserts these observational measurements
of the reference dynamical system into an approximate solution 
by means of an impulsive forcing.
In this way the approximating solution is coupled to the reference
solution at a discrete sequence of points in time.
This paper studies discrete data assimilation for the
Lorenz equations and the incompressible two-dimensional
Navier--Stokes equations.
In both cases we obtain bounds on the time interval $h$ between 
subsequent observations which guarantee the convergence of 
the approximating solution
obtained by discrete data assimilation to the reference solution.
\medskip
}%

\beginsection{Introduction}
In [\OLSTIT] and [\OLSTITB] Olson and Titi studied the number of 
determining modes for continuous data assimilation for the
incompressible two-dimensional Navier--Stokes equations.  
As in
those papers, the motivating problem for our work is the initialization
of weather forecasting models using near continuous in time measurement 
data obtained, for example, from satellite imaging.
In this work, rather than making the idealization that the measurement
data is continuous in time, we focus on the case where 
the measurement data is taken at a sequence
of discrete times $t_n$.

If $t_n=hn$ and $h$ is small then discrete data assimilation can be
viewed as near continuous.  One expects that if the approximating 
solution obtained by continuous data assimilation converges to the
reference solution then the approximating solution obtained by discrete
data assimilation for small $h$ will also converge to the reference
solution.  Note, however, that near continuous observational data
is mathematically quite different from continuous data.
If the observations are known continuously in time on some interval then
mathematically the $n$-th time derivatives may be calculated on that same
interval for all values of $n$.  It is possible that this derivative
information could lead to a reconstruction of the reference solution
in cases where near continuous measurement information might not.
For example, Wingard~[\WINGARD] shows for the Lorenz equations that
knowing $X$ and its time derivatives at a single point in time can be
used to recover both $Y$ and $Z$.

In this paper we present a technique to prove a discrete in time
determining mode result in the specific context of creating an
approximating solution that asymptotically converges to a reference
solution.  In addition to providing a more realistic framework
in which to study near continuous data assimilation, our work
can be seen as a discrete in time extension of the theory of
determining modes developed by Foias and Prodi [\FOIPRO] and
further refined by Jones and Titi in [\JONTITA] and [\JONTITB].
Unless otherwise noted, we shall assume the dynamics 
governing the evolution of the reference solution admit a global
attractor and that the reference
solution lies on that global attractor.
This assumption is made for simplicity, as our analysis 
actually depends only on the existence of an absorbing ball
which contains the reference solution forward in time.

We begin our discussion with the easier case
of the Lorenz equations to provide insight and illustrate the 
methods we will use for the incompressible two-dimensional
Navier--Stokes equations.  
General studies for the synchronization
of discrete in time coupled systems for the Lorenz equations were
produced by Yang, Yang and Yang [\YANG] and Wu, Lu, Wang and Liu [\WU].
Note, however, that the matrix corresponding to the observational
measurements studied here does not have a suitable spectral radius to 
apply the theorems of their work.

The method of discrete data assimilation can be described mathematically
as follows.
Let $U$ be a solution lying on the global attractor of a 
dissipative continuous dynamical system with initial condition 
$U_0$ at time $t_0$.
Let $S$ be the 
continuous semigroup defined by $U(t)=S(t,t_0,U_0)$.
Represent the observational 
measurements of the reference solution $U$ at time $t_n$ by $PU(t_n)$,
where $P$ is a finite-rank orthogonal projection and $t_n$ is an increasing 
sequence in time.
Discrete data assimilation inserts the observational measurements
into an approximate solution $u$ as the approximate solution 
is integrated in time.  
In particular, 
let $u_0=\eta+P U(t_0)$ and $u_{n+1}=QS(t_{n+1},t_n,u_n)+PU(t_{n+1})$ 
for $n=0,1,\ldots,$ where
$P\eta=0$ and $Q=I-P$.
Here $\eta$ corresponds to an initial guess for the part of 
the reference solution $QU(t_0)$ that can not be measured.
The approximating solution $u$ obtained by discrete data assimilation
is defined to be the piecewise continuous in time function 
$$u(t)=S(t,t_n,u_n)\quad\hbox{for}\quad t\in [t_n,t_{n+1}).
\eqno(\eqn{DASIM})$$
Our goal is to find conditions on $P$, $t_n$ and $\eta$ which guarantee 
that the approximating solution $u$
converges to the reference solution $U$ as $t\to\infty$.

For the Lorenz system the reference solution is a 
three dimensional vector consisting of the components
$X$, $Y$ and $Z$
whose evolution is governed by the coupled system of three ordinary
differential equations
$$\left\{\eqalign{
	\dot X&=-\sigma X+\sigma Y\cr
	\dot Y&=-\sigma X-Y-XZ\cr
	\dot Z&=-bZ+XY-b(r+\sigma)\cr
}\right. \eqno(\eqn{LORINTRO})
$$
where $\sigma=10$, $b=8/3$ and $r=28$.  We shall assume that the 
reference solution lies on the global attractor.

We take the observational measurements of the reference solution
to be the values of the variable $X$ at the times $t_n$.
These observations of $X$ are used to create an approximating 
solution whose components are $x$, $y$ and $z$.  Note that
$x(t_n)=X(t_n)$ where $y(t)$ and $z(t)$ are continuous
at $t=t_n$ and $x$, $y$ and $z$ satisfy
$$\left\{\eqalign{
	\dot x&=-\sigma x+\sigma y\cr
	\dot y&=-\sigma x-y-xz\cr
	\dot z&=-bz+xy-b(r+\sigma)\cr
}\right. \eqno(\eqn{LORINTROS})
$$
on each interval $[t_n,t_{n-1})$ for $n=0,1,2,\ldots\,$.

For simplicity we take $t_n=hn$ for some fixed $h>0$.  
Numerical experiments for the Lorenz system done by 
Hayden~[\KEVIN] indicate that the approximating solution
converges to the reference solution 
as $t\to\infty$ for values of $h$ as large as $0.175$.
\XDISCS\ shows the convergence of the approximating solution
to the reference solution when $h=0.1$.
\par
\noindent
\vbox{
\centerline{\epsfxsize=5truein\epsfbox{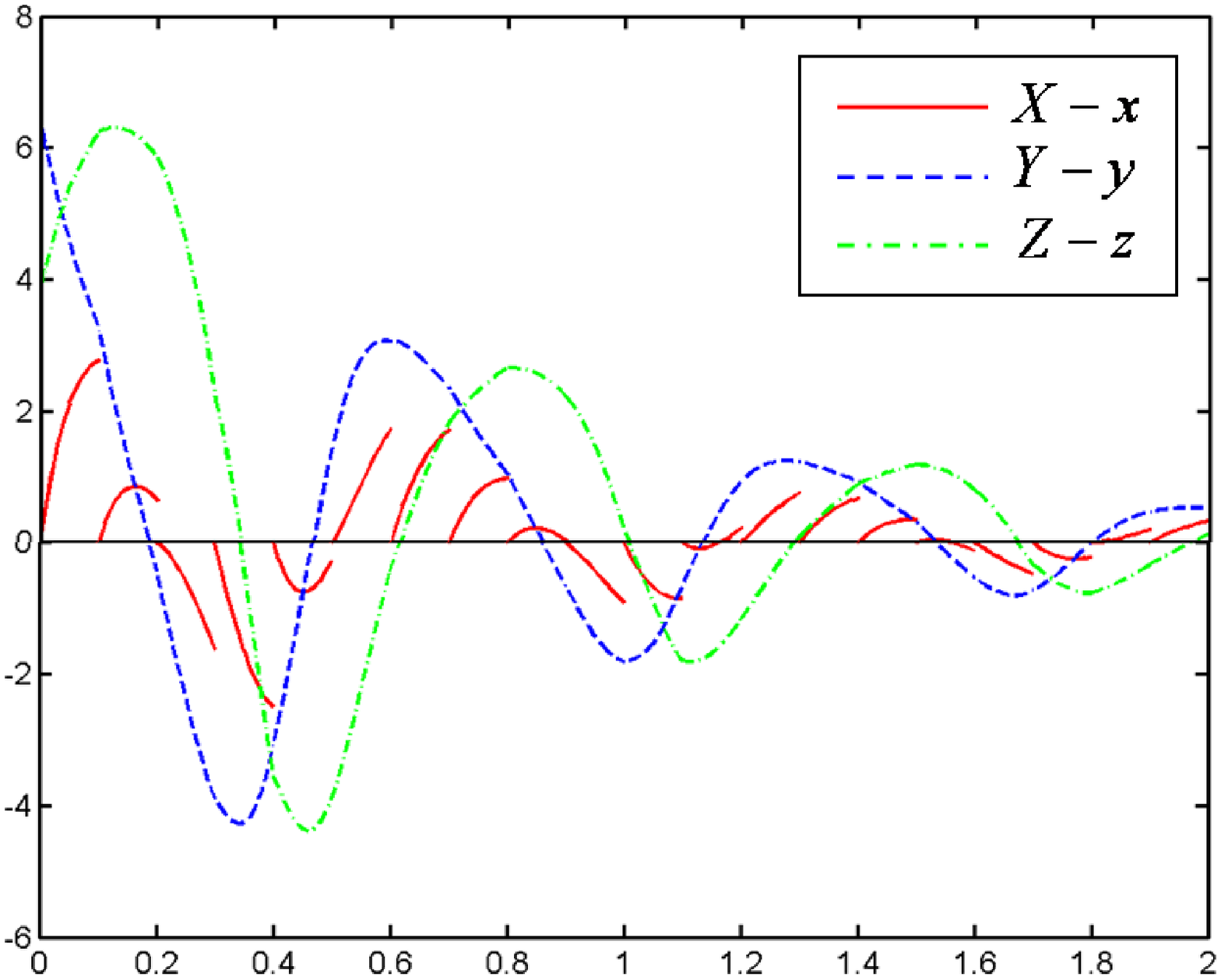}}
\narrower\narrower
\fig{XDISCS} Convergence of the approximating solution
to the reference solution for the Lorenz system when $h=0.1$.\par
}
\bigskip\noindent
Note again that the variables $y$ and $z$ in the 
approximating solution are continuous in time, in particular
at $t=t_n$, whereas $x$
is discontinuous.  The discontinuities in $x$ are the result of
the assimilation of the observations of $X$ at each time $t_n$.

Our main result for the Lorenz equations is an analytic proof of
\LORENZDA\ which shows there exists $t^*>0$ depending on $\sigma$,
$\beta$ and $r$ such that 
for any $h\in (0,t^*]$ the approximating solution obtained by discrete data 
assimilation of measurements of the $X$ variable at times $t_n=hn$,
as described in (\LORINTROS) above,
converges to the reference solution as $t\to\infty$.  

The second part of this paper focuses on the incompressible 
two-dimensional Navier--Stokes equations with $L$-periodic boundary
conditions
$$
    \left\{\eqalign{ {\partial U\over \partial t}
        -\nu\Delta U+(U\cdot \nabla)U+\nabla P=f\cr
        \nabla\cdot U=0\cr}\right.
\eqno(\eqn{TWODNS})
$$
where $\nu$ is the kinematic viscosity and $f$ is a time independent 
body forcing.
In this case $U$ can be expressed in terms of the 
Fourier series
$$
	U=\sum_{k\in\J} U_k e^{ik\cdot x} \quad\hbox{where}\quad
	\J=\left\{ {2\pi\over L}(n_1,n_2) :
		n_i\in\Z\hbox{ and } (n_1,n_2)\ne (0,0)\right\}
$$
where $U_k\cdot k=0$ for all $k\in\J$ and $U_k=U_{-k}^*$.
Again, we assume that the reference solution lies on the global attractor.

Let $P_\lambda$ be the orthogonal projection defined by
$$
	P_\lambda U = \sum_{k^2\le\lambda} U_k e^{ik\cdot x}.
	\eqno(\eqn{PLAMBDAD})
$$
We take the observational measurements of the reference solution
to be the values of $P_\lambda U$ at the times $t_n$.
Note that $\lambda^{-1/2}$ represents the smallest length scale of the 
fluid which can be observed---the resolution of the presumed measuring
equipment or distance between measuring stations.

In the context of continuous data assimilation and determining modes,
the rank of the smallest projection $P_\lambda$
such that the approximating solution converges to the reference
solution is called the number of determining modes.  For discrete
and near continuous data assimilation, the parameter $\lambda$
also depends on the interval of time between the observational measurements.  
Our main results for the incompressible two-dimensional Navier--Stokes 
equations are \BIGH\ and \ZEROETA.

\BIGH\ shows for any time interval $h>0$ there exists $\lambda$ large 
enough so that the approximating solution obtained by discrete data 
assimilation of the measurements $P_\lambda U(t_n)$ where $t_n=hn$ 
will converge to the reference solution.  This means that increasing
the resolution of the measurements can compensate for
a large time interval between subsequent observations.
\ZEROETA\ shows there is
a dimensionless constant $C$ such that if 
$$
 \lambda>{C\over \lambda_1^{5/3}}\Big({\|f\|_{L^2}\over \nu^2}\Big)^{8/3}
$$
then there exists $h>0$ depending only on $|f|$, $\nu$, $\Omega$
and $\lambda$ small enough such that
the approximating
solution obtained by discrete data assimilation with
initial guess $\eta=0$ 
converges to the reference solution.
This means that decreasing the time interval between subsequent observations
can compensate for a low resolution
provided the resolution meets a minimum standard.


\beginsection{Lorenz Equations}
Following Foias, Jolly, Kukavica and Titi [\FJKT]
we write the Lorenz system (\LORINTRO) as
$$
	{dU\over dt}+AU+B(U,U)=f \eqno(\eqn{LORENZ})
$$
where
$$
	U=\left[\matrix{X\cr Y\cr Z}\right],\quad
	A=\left[\matrix{ \sigma&-\sigma&0\cr
				\sigma&1&0\cr
				0&0&b\cr}\right],\quad
	f=\left[\matrix{0\cr 0\cr -b(r+\sigma)}\right],
$$
$$
	B\left(
		\left[\matrix{X\cr Y\cr Z}\right],
		\left[\matrix{\tilde X\cr \tilde Y\cr \tilde Z}\right]\right)
		=\left[\matrix{0\cr \phantom{-}(X\tilde Z+Z \tilde X)/2\cr 
			-(X\tilde Y+Y\tilde X)/2}\right]
$$
and again $\sigma=10$, $b=8/3$ and $r=28$.  One reason for
writing the Lorenz equations in this way is to make the similarities
and differences in the proofs from this section and following section
on the Navier--Stokes equations more transparent.
Despite the notational similarities there should be no trouble
distinguishing the results on the Lorenz equation that
apply only to the Lorenz equations from the results 
on the Navier--Stokes equations which apply only to the 
Navier--Stokes equations.  We start with some definitions and facts
that are easy to verify.

\dfn{NORMT} $|U|=\sqrt{(U,U)}=\sqrt{X^2+Y^2+Z^2}$.
\edfn

The facts below may be deduced from the preceding definitions and 
are listed here for reference.
First, we state an estimate involving the linear term $AU$ that plays the
same role in our treatment of the Lorenz system that \POINCARE\ plays
in our later treatment of the incompressible two-dimensional Navier--Stokes
equations: 
$$(AU,U)=\sigma X^2+Y^2+bZ^2\ge |U|^2. \eqno(\eqn{POINL})$$ 
Next we state some algebraic identities analogous to the orthogonality
relations (\ORTHAA) and (\ORTHAB) for the nonlinear term in the 
two-dimensional Navier--Stokes equations:
$$(B(U,U),U)=0 \qquad\hbox{and}\qquad
B(U,\tilde U)=B(\tilde U,U). \eqno(\eqn{SYMML}) $$
We will also use an estimate on $B$ which is similar to \NONLIN:
$$|B(U,\tilde U)|\le 2^{-1} |U||\tilde U|. \eqno(\eqn{NONLINL})$$

The following bound on the global attractor was shown in
Temam [\TEMA] on page 33, see also Foias, Constantin and Temam [\CONFOITEM].
\thm{BOUNDK}  Let $U$ be a trajectory that lies on the global 
attractor of (\LORENZ).  Then
$|U(t)|^2\le K$ for all $t\in\R$ where 
$$K={b^2(r+\sigma)^2\over 4(b-1)}.
\eqno(\eqn{ETWO})
$$
\ethm

Before proceeding, we state Young's inequality which will be used throughout
the remainder of this work.
\thm{YOUNGS} Let 
${1/p}+{1/q}=1$
then $|xy|\le |x|^p/p+|y|^q/q.$
\ethm

We are now ready to begin our study of discrete data assimilation
for the Lorenz equations.
Our measurements will consist of the $X$ variable
at the times $t_n$.  Therefore, we define the orthogonal
projection $P$ as
$$P=\left[\matrix{1&0&0\cr 0&0&0\cr 0&0&0}\right]\qquad
\hbox{and}\qquad
Q=I-P=\left[\matrix{0&0&0\cr 0&1&0\cr 0&0&1}\right].$$

Let $U$ be a solution of the Lorenz 
equations lying on the global attractor and $u$ the approximating 
solution given by (\DASIM) where $S$ is the semigroup 
generated by (\LORENZ).
Since $U$ and $u$ both satisfy (\LORENZ) on the time
interval $[t_n,t_{n+1})$ 
setting $\delta=U-u$ we obtain
$$
	{d\delta\over dt}+\nu A\delta +B(U,U)-B(u,u)=0
\quad\hbox{for}\quad t\in[t_n,t_{n+1})$$
or after some algebra that
$$
	{d\delta\over dt}+ A\delta +B(U,\delta)
		+B(\delta,U)-B(\delta,\delta)=0
\quad\hbox{for}\quad t\in[t_n,t_{n+1}).
\eqno(\eqn{DIFEQL})
$$

\lem{EXGROWTHL} 
There exists $\beta>0$ given by (\EONE) below, depending on $\sigma$,
$b$ and $r$, such that 
$|\delta(t)|^2\le |\delta(t_n)|^2 e^{\beta (t-t_n)}$
for $t\in [t_n,t_{n+1})$.
\prf
Take the inner product of $\delta$ with (\DIFEQL) and apply
(\POINL) and (\SYMML) to obtain
$$
	{1\over 2} {d |\delta|^2\over dt}
		+|\delta|^2+2(B(U,\delta),\delta)\le 0.
\eqno(\eqn{EXGA})
$$
Estimating using (\NONLINL) and \BOUNDK\ gives
$$
	2|(B(U,\delta),\delta)|\le |U||\delta|^2
		\le K^{1/2} |\delta|^2.
$$
Therefore
$$
	{d |\delta|^2\over dt} \le \beta |\delta|^2
$$
where 
$$\beta=2 (K^{1/2}-1)\eqno(\eqn{EONE})$$
and $K$ is defined in (\ETWO).  Integrating
from $t_n$ to $t$ yields
$$
	|\delta(t)|^2\le |\delta(t_n)|^2 e^{\beta (t-t_n)},
$$
which finishes the proof.
\eprfn

\thm{LORENZDA}
Let $U$ be a solution of the the Lorenz equations (\LORENZ)
lying on the global attractor.
Then, there exists $t^*>0$ depending 
only on $\sigma$, $b$ and $r$
such that for any $h\in(0,t^*]$ the approximating solution $u$ given
by (\DASIM), see also (\LORINTROS), with $t_n=hn$ converges 
to $U$ as $t\to\infty$.
\prf
Take the inner product of $P\delta$ with (\DIFEQL) to obtain
$$
	{1\over 2}
	{d |P\delta|^2\over dt}+ (A\delta,P\delta)=0.
\eqno(\eqn{LORDE})
$$
Let
$$
	w=(PA-AP) \delta= \left[\matrix{0&-\sigma&0\cr
		-\sigma&0&0\cr 0&0&0}\right] \delta
		=-\sigma\left[\matrix{\delta_Y\cr\delta_X\cr0}\right]
\quad\hbox{where}\quad
	\delta=\left[\matrix{\delta_X\cr\delta_Y\cr\delta_Z}\right].
$$
Then $|P\delta|^2=\delta_X^2$ and
$$
	(A\delta,P\delta)=
		(PA\delta,P\delta)=
			(AP\delta,P\delta)+(w,P\delta)
			=\sigma \delta_X^2-\sigma \delta_X\delta_Y.
$$	
Therefore (\LORDE) becomes
$$
	{1\over 2}
	{d \delta_X^2\over dt}+\sigma\delta_X^2=\sigma \delta_X\delta_Y
		\le {\sigma\over 2} \delta_X^2+{\sigma\over 2}\delta_Y^2
		\le {\sigma\over 2} \delta_X^2+{\sigma\over 2}|\delta|^2
$$
where we have applied \YOUNGS\ with $p=q=2$.
Using \EXGROWTHL\ we obtain
$$
	{d |P\delta|^2\over dt}+\sigma|P\delta|^2
		\le \sigma \delta^2
		\le \sigma |\delta(t_n)|^2 e^{\beta(t-t_n)}.
$$
Multiplying by the integrating factor $e^{\sigma (t-t_n)}$ and
integrating from $t_n$ to $t$ gives
$$
	|P\delta(t)|^2\le
		{\sigma |\delta(t_n)|^2\over \beta+\sigma}
		\big(
			e^{\beta(t-t_n)}-e^{-\sigma(t-t_n)}
		\big)
\eqno(\eqn{THISB})
$$
where we have used the fact that $P\delta(t_n)=0$.

A finer analysis of the nonlinear term appearing in (\EXGA) gives
$$
	2(B(U,\delta),\delta)=\delta_X(Z\delta_Y -Y\delta_Z).
$$
Therefore
$$|2(B(U,\delta),\delta)|
	\le |U||P\delta||\delta|\le K^{1/2} |P\delta||\delta|
	\le {1\over 2}K |P\delta|^2+ {1\over 2}|\delta|^2.
\eqno(\eqn{THISA})
$$
Substituting (\THISB) and (\THISA) into (\EXGA) yields
$$
	{d |\delta|^2\over dt}
		+|\delta|^2\le K 
		{\sigma |\delta(t_n)|^2\over \beta+\sigma}
		\big(
			e^{\beta(t-t_n)}-e^{-\sigma(t-t_n)}
		\big).
	\eqno(\eqn{THISC})
$$
Multiply (\THISC) by $e^{(t-t_n)}$ and
integrate from $t_n$ to $t$ to obtain
$$
	|\delta(t)|^2
		\le M(t-t_n)|\delta(t_n)|^2
$$
where
$$
	M(\tau)=
	 e^{-\tau} 
	\Big(1+{\sigma K\over \beta+\sigma}
		\int_{0}^\tau
		\big(
            e^{(\beta+1)s}-e^{-(\sigma-1)s}
		\big) ds\Big)
	\eqno(\eqn{THISD})
$$
is a function that doesn't depend on $t_n$.
Note that $M(0)=1$.  Differentiating yields
$$
	M'(\tau)=
	-M(\tau) + 
	{\sigma K\over \beta+\sigma}
		\big(
            e^{\beta\tau}-e^{-\sigma\tau}
		\big).
$$
Therefore $M'(0)=-1$.  It follows that there is $t^*>0$ 
such that $$M(h)=1+\int_0^{h} M'(s) ds<1$$
for all $h\in (0,t^*]$.

Next fix $h\in (0,t^*]$ and let $\gamma=M(h)<1$.  Then
$$
\eqalign{
	|\delta(t_{n+1})|^2
		&=|Q\delta(t_{n+1})|^2
		=\lim_{t\nearrow t_{n+1}} |Q\delta(t)|^2
		\le \lim_{t\nearrow t_{n+1}} |\delta(t)|^2\cr
		&\le \lim_{t\nearrow t_{n+1}} M(t-t_{n})|\delta(t_{n})|^2
		= M(h)|\delta(t_{n})|^2=\gamma |\delta(t_{n})|^2
}$$
implies by induction that
$$ |\delta(t_{n})|^2
	\le \gamma^n |\delta(t_0)|^2
	=\gamma^n|QU(t_0)-\eta|^2
	\le \gamma^n R
$$
where $R= 2 (K+|\eta|^2).$
Now let $t>0$ and choose $n$ so that $t\in [t_n,t_{n+1})$.
Then $n\to\infty$ as $t\to\infty$ and therefore
$$|\delta(t)|^2\le M(t-t_n)|\delta(t_n)|^2
		\le |\delta(t_n)|^2\le \gamma^n R\to 0
\eqno(\eqn{UBOUNDB})
$$
shows that the approximating solution converges to 
the reference solution as $t\to\infty$.
\eprfn

\cor{ANBND}
If $\sigma=10$, $b=8/3$ and $r=28$, then
$t^*\approx 0.000129$.
\prf
The result follows using the value of $K$ from \BOUNDK\ 
and choosing $t^*$ slightly less than the value of $t$ such 
that $M(t)=1$.
\eprf

Before proving the final result in this section, we extend
\LORENZDA\ by proving 
that the interval of time between $t_{n+1}$ and $t_n$ 
need not be exactly $h$ for every $n\in \N$.

\cor{VARHL}  Let $t^*$ be the bound given in \LORENZDA.  Suppose
$t_{n+1}-t_{n}\le t^*$ where $t_n\to\infty$ as $n\to\infty$.  Then the 
approximating solution $u$ given by (\DASIM) converges to the reference
solution $U$ of the Lorenz equations (\LORENZ) as $t\to\infty$.
\prf
Let $h_n=t_{n}-t_{n-1}$. If there exists $\epsilon>0$
such that the set ${\cal K}=\{\,k: h_k\ge \epsilon\,\}$ 
is infinite then
$M(h_k)\le \max\big\{ \,M(s): s\in [\epsilon,t^*]\, \big\}<1$
for $k\in{\cal K}$ implies that
$$
|\delta(t_n)|^2\le \prod_{k=1}^{n} M(h_k) R\to 0
\quad\hbox{as}\quad n\to\infty.
\eqno(\eqn{DELTACON})
$$
Otherwise, $h_n\to 0$ as $n\to\infty$.
By Taylor's theorem
$|M(h_n)-M(0)-h_nM'(0)|\le C h_n^2 $
where $C={1\over 2}\max\big\{\, |M''(s)|: s\in [0,t^*]\,\big\}.$
Choose $N$ so large that $Ch_k-1\le -1/2$ for $k\ge N$.
Since $M(0)=1$, $M'(0)=-1$ and $\log x\le x-1$ for $x>0$ it follows that
$$\eqalign{
	\sum_{k=N}^n\log M(h_k)
	& 
	\le \sum_{k=N}^n\big( M(h_k)-1\big)
	\le \sum_{k=N}^n \big( M(0)+h_k M'(0)+C h_k^2-1\big)\cr
	&= \sum_{k=N}^n h_k (C h_k -1)
	\le
    -{1\over 2}\sum_{k=N}^n h_k\to -\infty
	\quad\hbox{as}\quad n\to\infty.
}$$
Thus 
$$
	\prod_{k=1}^n M(h_k) R 
		=\exp\Big(\log R+
			\sum_{k=1}^n \log M(h_k)
		\Big)\to 0
\quad\hbox{as}\quad n\to\infty.
$$
The proof now finishes as in \LORENZDA.
\eprf

Our final result on the Lorenz equations shows that the approximating
solution is bounded for any updating time interval of $h$.
If $h\le t^*$ then the approximating solution converges
to the reference solution, and since the reference solution
is bounded, then the approximating solution will also be bounded.
In the case where the approximating solution
does not converge to the reference solution then the following
result shows that the approximating solution is still bounded.

\thm{BNDSLAVE} Let $U$ be a trajectory that lies on the
global attractor of~(\LORENZ).  The 
approximating solution $u$ given by (\DASIM) where
$t_n=nh$ with $h>0$ 
is bounded.  Namely,
there is  a constant $M_1$ that depends
only on $\eta$, $\sigma$, $b$ and $r$ such that
$|u(t)|^2\le M_1/(1-e^{-h})$ for all $t\ge0$ and $h>0$.
\prf
Taking inner product of (\LORENZ) with $u$ and
using (\POINL) followed by Young's inequality we obtain
$$
	{1\over 2}{d |u|^2\over dt}+ |u|^2\le (f,u)
		\le {1\over 2}|f|^2+{1\over 2}|u|^2
$$
and consequently
$$
	{d |u|^2\over dt}+|u|^2\le |f|^2
	\eqno(\eqn{CANIMPROVE})
$$
for $t\in [t_n,t_{n+1})$.  Gr\"onwall's inequality then implies
$$
	|u(t)|^2\le |u_n|^2 e^{-(t-t_n)}+|f|^2(1-e^{-(t-t_n)})
\qquad\hbox{for}\qquad t\in [t_n,t_{n+1}).$$

Defining
$$
	\bar u_0=\eta
\qquad\hbox{and}\qquad
	\bar u_{n+1}=
\lim_{t\nearrow t_{n+1}} u(t)= S(t_{n+1},t_n,u_n)$$
so that $u_n=Q\bar u_n+ PU(t_n)$ we obtain
$$\eqalign{
	|\bar u_{n+1}|^2
    &\le \lim_{t\nearrow t_{n+1}} \big(
	|u_n|^2 e^{-(t-t_n)}+|f|^2(1-e^{-(t-t_n)})\big)\cr
    &\le |u_n|^2 \gamma+|f|^2(1-\gamma)}
$$
where $\gamma=e^{-h}$.
Since
$$\eqalign{
	|u_n|^2=|Q\bar u_n|^2+|PU(t_n)|^2
		\le|\bar u_n|^2+|U(t_n)|^2\le |\bar u_n|^2+K
}$$
we obtain 
$$|\bar u_{n+1}|^2 \le 
	|\bar u_n|^2\gamma +C_1
\eqno(\eqn{INDUL})
$$
where $C_1=K\gamma +|f|^2(1-\gamma).$

Induction on (\INDUL) and summing the series yields
$$
	|\bar u_n|^2\le  |\bar u_0|^2\gamma^n
		+ C_1(1+\gamma+\cdots+\gamma^{n-1})
		= |\eta|^2\gamma^n+C_1{1-\gamma^n\over 1-\gamma}.$$
Given $t>0$ choose $n$ so that $t\in [t_n,t_{n+1})$.  Then
$t=t_n+\alpha$ where $\alpha\in [0,h)$ and
$$\eqalign{
	|u(t)|^2
	&\le |u_n|^2 e^{-\alpha}+|f|^2(1-e^{-\alpha})
	\le |u_n|^2+|f|^2\cr
	&\le |\bar u_n|^2 + K + |f|^2
	\le  |\eta|^2 \gamma^n+ C_1 {1-\gamma^n\over 1-\gamma}+ 
		K + |f|^2
	\le {M_1\over 1-\gamma}}$$
where
$
	M_1=|\eta|^2+C_1+K+|f|^2.$
\eprf

If we take $\eta=0$ as our initial guess for
$QU(t_0)$ when forming the approximating solution $u$ then 
the constant $M_1$ depends only on $\sigma$, $b$ and $r$.  In
either case we obtain the asymptotic bound
$$
	\limsup_{t\to\infty} |u(t)|^2 \le {M_2\over 1-e^{-h}}.
$$
where the constant $M_2$ depends only on $\sigma$, $b$ and $r$.

\BNDSLAVE\ can be improved using 
the exact structure of $(Au,u)$ and $(f,u)$.
In particular, we obtain using \YOUNGS\ that
$$\eqalign{
    (Au,u)-(f,u)
	&= \sigma x^2+y^2+bz^2+ zb(r+\sigma)\cr
	&= \sigma x^2+y^2+z^2+(b-1)z^2+zb(r+\sigma)\cr
        &\ge \sigma x^2+y^2+z^2-{b^2(r+\sigma)^2\over 4(b-1)}
        \ge |u|^2-{|f|^2\over 4(b-1)}}
$$
which improves the bound to $|u(t)|^2\le M_3/(1-e^{-2h})$.
Although an improvement, this bound on $u$ still tends to infinity
as $h$ tends to zero.  As mentioned earlier, a bound 
uniform in $h$ can be obtained
by combining \LORENZDA\ with \BNDSLAVE.

\cor{UBOUND}  There exists a bound $M_4$ depending only on $\eta$,
$\sigma$, $b$ and $r$ such that the approximate solution $u$ obtained
with $t_n=hn$ is bounded by $M_4$ for any $h>0$.
\prf
Let $t^*$ be given as in \LORENZDA.  If $h<t^*$ then (\UBOUNDB)
implies
$$
	|u(t)|=|U(t)-\delta(t)|\le |U(t)|+|\delta(t)|
		\le K^{1/2}+R^{1/2}.
$$
If $h\ge t^*$ then \BNDSLAVE\ implies
$$
	|u(t)|\le M_1^{1/2}/(1-e^{-h})^{1/2}
		\le M_1^{1/2}/(1-e^{-t^*})^{1/2}.
$$
Taking 
	$$M_4=\max\{ K^{1/2}+R^{1/2}, M_1^{1/2}/(1-e^{-t^*})^{1/2}\}$$
yields a constant that depends only on $\eta$, $\sigma$, $b$ and $r$.
\eprf

As mentioned before, if we 
take $\eta=0$ or consider asymptotic bounds
on $u$ as $t\to\infty$ then the dependency on $\eta$ can be
removed from these bounds.

\beginsection{Navier--Stokes Equations}
This section contains results for the two-dimensional
incompressible Navier--Stokes equations that are similar
to the results proved in the previous section for
the Lorenz equations.
Let 
$\Omega=[0,L]\times [0,L]$,
${\cal V}$ be the space of divergence-free vector-valued 
$L$-periodic trigonometric polynomials from $\Omega$
into $\R^2$ with $\int_\Omega u=0$, 
$H$ be the closure of ${\cal V}$
with respect to the $L^2$ norm and
$P_H$ be the $L^2$-orthogonal projection $P_H\colon L^2(\Omega)\to H$,
referred to as the Leray-Helmholtz projector.
Following the notations of Constantin and Foias [\CONFOI],
Robinson~[\ROB] and Temam [\TEMB] we write the incompressible 
two-dimensional Navier--Stokes equations (\TWODNS) as
$$
	{du\over dt} +\nu Au+B(u,u)=f\eqno(\eqn{FTWODNS})
$$
where
$$
	Au=-P_H\Delta u\qquad\hbox{and}\qquad
	B(u,v)=P_H\big[(u\cdot \nabla)v\big].
	\eqno(\eqn{DEFAB})
$$
Note that we have assumed $f\in H$ so that $P_H f=f$.
Also notice that $A=-\Delta$ in the periodic case.

Let 
$V$ be the closure of $\V$ with respect to the $H^1$ norm and
${\cal D}(A)=H\cap H^2(\Omega)$ be the domain of $A$.
We may define norms on $H$, $V$ and ${\cal D}(A)$ which are
equivalent to the $L^2$, $H^1$ and $H^2$ norms respectively by
$$
|u|=L^2\sum_{k\in{\cal J}} |u_k|^2, \qquad
\|u\|=L^2\sum_{k\in{\cal J}}k^2 |u_k|^2 \qquad \hbox{and} \qquad
|Au|=L^2\sum_{k\in{\cal J}}k^4 |u_k|^2.
$$ 
Here $u$ has been expressed in terms of the Fourier series
$$
	u=\sum_{k\in\J} u_k e^{ik\cdot x} \quad\hbox{where}\quad
	\J=\left\{ {2\pi\over L}(n_1,n_2) :
		n_i\in\Z\hbox{ and } (n_1,n_2)\ne 0\right\}.
$$

The mathematical theory proving the existence and uniqueness
of strong solutions to the two-dimensional incompressible
Navier--Stokes equations (\FTWODNS) with initial data in $V$
may be found for example in
[\CONFOI], [\ROB] or [\TEMB].  Specifically we have

\thm{STRONGSOL} Let $u_0\in V$
and $f\in H$.  Then (\FTWODNS) has unique strong
solutions that satisfy
$$
	u\in L^{\infty}\big((0,T);V\big)
		\cap L^2\big((0,T);{\cal D}(A)\big)
\quad\hbox{and}\quad
	{du\over dt} \in L^2\big((0,T);H\big)
$$
for any $T>0$.  Furthermore, this solution is in
$C\big([0,T];V\big)$ and depends continuously
on the initial data $u_0$ in the $V$ norm.
\ethm

Let $S(t,t_0,u_0)$ for $t\ge t_0$ denote the unique strong solution 
of (\FTWODNS) given by \STRONGSOL\ 
with initial condition $u_0$ at time $t_0$.
Let
$\lambda_1=(2\pi)^2/L^2$ be the smallest eigenvalue of $A$ on 
${\cal D}(A)$ and define
$P_\lambda$ as in~(\PLAMBDAD).
Further define $Q_\lambda=I-P_\lambda$ to be the orthogonal
complement of $P_\lambda$.
The Poincar\'e inequalities for $u$,
$P_\lambda u$ and $Q_\lambda u$ on $\Omega$ can now be
summarized as

\thm{POINCARE} Given $P_\lambda$, $Q_\lambda$ and $\lambda_1$
as defined above then 
$$
		|u|\le \lambda_1^{-1/2} \|u\|,\qquad
		|Q_\lambda u|\le \lambda^{-1/2} \|Q_\lambda u\|,\qquad
		\|P_\lambda u\|\le \lambda^{1/2} |P_\lambda u|$$
and
$$
		\|u\|\le \lambda_1^{-1/2} |Au|,\qquad
		\|Q_\lambda u\|\le \lambda^{-1/2} |AQ_\lambda u|,\qquad
		|AP_\lambda u\|\le \lambda^{1/2} \|P_\lambda u\|.$$
provided the norms exist and are finite.


\thm{PYTHAG} If $u\in {\cal D}(A)$ then
$$ |Au|^{1/4}\le 2^{1/8}(\lambda^{1/8} \|u\|^{1/4}+|Q_\lambda Au|^{1/4}).$$
\prf
Since $P_\lambda$ and $Q_\lambda$ are orthogonal then
$$
	|Au|^2=|P_\lambda A u|^2+|Q_\lambda A u|^2
		\le \lambda \|P_\lambda u\|^2+|Q_\lambda Au|^2
	\le 2\max\{\lambda \|P_\lambda u\|^2,|Q_\lambda A u|^2\}.
$$
Therefore
$$
	|Au|^{1/4}\le 2^{1/8}\max\{\lambda^{1/8} \|P_\lambda u\|^{1/4}
			,|Q_\lambda A u|^{1/4}\}
		\le  2^{1/8}(\lambda^{1/8} \|u\|^{1/4}+|Q_\lambda Au|^{1/4}),
$$
which completes the proof.
\eprf

Let us now recall some algebraic properties 
of the non-linear term that can also be found 
in [\CONFOI], [\ROB] or~[\TEMB]
that play an important role in our analysis.  They are
$$
	(B(u,v),w)= (B(u,w),v), \eqno(\eqn{ORTHAA})
$$
$$
	(B(u,v),v)=0 \qquad \hbox{and} \qquad
	(B(u,u),Au)=0.
	\eqno(\eqn{ORTHAB})
$$

We now move on to some inequalities 
which we shall refer to later.

\thm{SOBOLEV}  There exists a dimensionless
constant $c$ such that if $u\in V$ then 
$$ 
	\|u\|_{L^{8/3}}^2\le c |u|^{3/2} \|u\|^{1/2}, \qquad
	\|u\|_{L^4}^2\le c |u|\|u\|, \qquad
	\|u\|_{L^8}^2\le c |u|^{1/2} \|u\|^{3/2}
$$
and if $u\in {\cal D(A)}$ then
$$
	\|u\|_{L^\infty}
		\le c \lambda_1^{-1/8} \|u\|^{3/4} |Au|^{1/4}.
$$

The first three inequalities above may be obtained from the Sobolev 
inequalities followed by interpolation, see for example (6.2) and
(6.7) in [\CONFOI].  The inequality bounding
$L^4$ is sometimes referred to as Ladyzhenskaya's inequality and
appears as Lemma 1 on page~8 of Ladyzhenskaya [\LADY].
The last inequality is a form of Agmon's inequality
which appears as (2.23) on page 11 in Temam [\TEMB].
These inequalities may be used along with H\"older's inequality
to estimate the nonlinear term.

\thm{NONLIN} Let $u,v,w\in V$ then
$$
	|(B(u,v),w)|\le c |u|^{1/2}\|u\|^{1/2} \|v\| |w|^{1/2} \|w\|^{1/2}
$$
if further $v\in {\cal D}(A)$ then
$$\eqalign{
	|(B(u,v),w)|&\le
		c\lambda_1^{-1/8} \|u\| \|v\|^{3/4} |Av|^{1/4} |w|\cr
	|(B(v,u),w)|&\le
		c\lambda_1^{-1/8} \|u\| \|v\|^{3/4} |Av|^{1/4} |w|}
$$
where $c$ is the constant appearing in \SOBOLEV.
\ethm

We consider a reference solution $U$ 
to the incompressible two-dimensional Navier--Stokes equations
and the approximating solution $u$ given by (\DASIM) where
$S$ is the semigroup generated by (\FTWODNS).  As with the
Lorenz equations, we shall assume that the reference solution
lies on the global attractor.  As proved in Jones and Titi [\JONTITA] we have

\thm{APRIORI} A solution $U$ that lies on the global attractor 
of (\FTWODNS) satisfies the bound  $\|U\|^2\le K$
where $$K={|f|^2\over \lambda_1\nu^2}.\eqno(\eqn{EFOUR})$$
\ethm

A similar result also appears in [\CONFOI] for establishing 
estimates on the global attractor.  As mentioned in
the introduction, any reference solution $U$ 
which satisfies a bound such as~(\EFOUR) forward in time
is suitable for our analysis, whether that solution is on 
the attractor or not.  However, for simplicity we continue
to assume $U$ lies on the global attractor.

Now setting $\delta=U-u$ and following the same algebra as 
for the Lorenz equations we arrive at 
$$
    {d\delta\over dt}+\nu A\delta +B(U,\delta)
        +B(\delta,U)-B(\delta,\delta)=0,
\eqno(\eqn{DIFEQ})
$$
an equation that looks like (\DIFEQL) where $A$ 
and $B$ have been given the meanings in (\DEFAB).

Before starting with the proof of our main result we begin with
a simpler theorem that finds values of $\lambda$ large 
enough for the approximating solution $u$ to start converging to 
a reference solution $U$ lying on the attractor but does 
not provide the uniform bound on~$h$ necessary to maintain
convergence as $t\to\infty$.

\thm{NOINDUCT}
If $\lambda>c^2 |f|^2/(\lambda_1\nu^4)$ there exists $t>t_n$
such that $|\delta(t)|<|\delta(t_n)|$.
\prf
Multiply (\DIFEQ) by $\delta$ and integrate over $\Omega$.
The first of the orthogonality relationships given in (\ORTHAB) yields
$$
	{1\over 2}{d|\delta|^2\over dt}+\nu \|\delta\|^2 
		=-(B(\delta,U),\delta).
\eqno(\eqn{DELTAEQ})
$$
Now applying \NONLIN\ followed by Young's and \APRIORI\ we obtain
$$
	{1\over 2}{d|\delta|^2\over dt}+\nu \|\delta\|^2 
		\le c |\delta|\|\delta\| \|U\|
		\le {\nu\over 2} \|\delta\|^2+ {c^2 K\over 2\nu} |\delta|^2
$$
or
$$
	{d|\delta|^2\over dt}\le 
		{c^2 K\over \nu} |\delta|^2
		-\nu \|\delta\|^2.
$$
Integrating from $t_n$ to $t$ where $t\in[t_n,t_{n+1})$ results in
$$
	|\delta(t)|^2-|\delta(t_n)|^2\le
		\int_{t_n}^t \Big(
        	{c^2 K\over \nu} |\delta(s)|^2
			-\nu \|\delta(s)\|^2 \Big) ds.
$$
Therefore
$$
	|\delta(t)|^2\le M_n(t-t_n)|\delta(t_n)|^2
$$
where
$$
	M_n(\tau)=1+{1\over |\delta(t_n)|^2}
		\int_0^{t-t_n} \Big(
			{c^2 K\over \nu} |\delta(t_n+s)|^2
			-\nu \|\delta(t_n+s)\|^2 \Big) ds.$$
Note $M_n(\tau)$ is a differentiable function with $M_n(0)=1$.

Differentiating yields
$$
	M_n'(\tau)= {c^2 K\over \nu}
			{|\delta(t_n+\tau)|^2\over |\delta(t_n)|^2}
-\nu{ \|\delta(t_n+\tau)\|^2\over |\delta(t_n)|^2}. $$
By \POINCARE\ and the hypothesis $\lambda>c^2 K/\nu^2$ we obtain
$$
	M_n'(0)=
		{c^2 K\over \nu} 
		-\nu{ \|\delta(t_n)\|^2\over |\delta(t_n)|^2}
	=
		{c^2 K\over \nu} 
		-\nu{ \|Q_\lambda \delta(t_n)\|^2\over |\delta(t_n)|^2}
	\le
		{c^2 K\over \nu} -\nu\lambda <0.
$$
Therefore there exists $\tau>0$ such that
$$
	M_n(\tau)=1+\int_{0}^\tau M'(s) ds <1.
$$
It follows that for $t=t_n+\tau$ that $|\delta(t)|<|\delta(t_n)|$.
\eprf

It should be pointed out that the value $t$ in the
above theorem depends on $M_n(\tau)$ and thus on
$\delta(s)$ for $s\ge t_n$.
Unfortunately, this provides us with no way of knowing whether 
the family of functions $M_n$ is equicontinuous or not.
This means
\NOINDUCT\ can not be used to provide a uniform bound $t^*$ on 
the discrete measurement time interval $h$ that ensures the 
approximate solution 
will converge to the reference solution as $t\to\infty$.
In \NSDA\ below we circumvent this issue at the expense
of a more stringent condition on the 
minimum size of $\lambda$.

Note that
the bound on $\lambda$ given
in \NOINDUCT\ above is the same as the bound 
given by Jones and Titi in [\JONTITA]
on the number of determining modes for the incompressible
two-dimensional Navier--Stokes equations.
This bound was later improved by Jones and Titi in [\JONTITB] using 
estimates on the time averages of the term $\|Au\|^2$.
However, in the case of discrete assimilation the
same technique did not achieve similar improvements.

We commence with the detailed analysis that allows us to prove
our main result for the incompressible two-dimensional 
Navier--Stokes equations.  The first theorem we
prove is an analog of \EXGROWTHL.

\thm{EXGROWTH} There exists $\beta>0$ depending on
$|f|$, $\Omega$ and $\nu$ such that the solution 
$\delta$ to
(\DIFEQ) satisfies $\|\delta(t)\|^2\le \|\delta(t_n)\|^2 e^{\beta (t-t_n)}$
for $t\in [t_n,t_{n+1})$.
\prf
Multiply (\DIFEQ) by $A\delta$ and integrate over $\Omega$ to
obtain
$$
	{1\over 2} {d\|\delta\|^2\over dt} + \nu|A\delta|^2
	+ (B(U,\delta),A\delta) +(B(\delta,U),A\delta)=0.$$
By \NONLIN\ we have that
$$
	|(B(\delta,U),A\delta)|\le
		c\lambda_1^{-1/8} \|\delta\|^{3/4}\|U\||A\delta|^{5/4}
$$
and
$$
	|(B(U,\delta),A\delta)|\le
		c\lambda_1^{-1/8} \|\delta\|^{3/4}\|U\||A\delta|^{5/4}.
$$
Therefore applying \YOUNGS\ with $p=8/3$ and $q=8/5$ we
obtain
$$\eqalign{
	{1\over 2} {d\|\delta\|^2\over dt} + \nu|A\delta|^2
	&\le
		2c\lambda_1^{-1/8} \|\delta\|^{3/4}\|U\||A\delta|^{5/4}\cr
	&\le
		C_1\nu^{-5/3}\lambda_1^{-1/3} 
			\|\delta\|^2\|U\|^{8/3}+ \nu|A\delta|^2\cr
	&\le
		C_1\nu^{-5/3}\lambda_1^{-1/3} K^{4/3}
			\|\delta\|^2 + \nu|A\delta|^2\cr}$$
where $C_1$ is a dimensionless constant related to $c$.  Hence
$$
	{d\|\delta\|^2\over dt}\le \beta \|\delta\|^2
$$
where $\beta=2C_1\nu^{-5/3}\lambda_1^{-1/3} K^{4/3}$.  Integrating
in time from $t_n$ to $t$ and applying Gr\"onwall's
inequality yields
$\|\delta(t)\|^2\le \|\delta(t_n)\|^2 e^{\beta (t-t_n)}$
for $t\in [t_n,t_{n+1})$.
\eprf

We are now ready to prove the main results of this paper.

\thm{NSDA}
If
$$	\lambda>{9\over \lambda_1^{1/3}}
		\Big({2c |f|/(\lambda^{1/2}\nu)+
		c\|\delta(t_0)\|\over\nu }\Big)^{8/3}
	\eqno(\eqn{NSDAL})
$$
then there exists $t^*>0$ depending only on
$K$, $\|\delta(t_0)\|$, $\nu$, $\Omega$ and 
$\lambda$ such that for any $h\in (0,t^*]$ the
approximating solution $u$ given by (\DASIM) with $t_n=hn$ converges 
to the reference solution $U$ of (\FTWODNS) as $t\to\infty$.
\prf
The proof is similar to the proof of \LORENZDA\ for the
Lorenz equations with the addition that 
we first use induction to show the bound $R=\|\delta(t_0)\|^2$ on 
the difference of the initial conditions ensures 
that $\|\delta(t_n)\|^2\le R$ holds for 
each each $t_n$.

Define
$$	
	g(\tau)=C_2\Big({\lambda\over \nu^4\lambda_1}\Big)^{1/4} g_1(\tau)
        +C_3\Big({1\over \nu^{5}\lambda_1}\Big)^{1/3} g_2(\tau)
	\eqno(\eqn{DEFGEE})
$$
where
$$
	C_2=c^2 (2^{5/4}), \qquad
	C_3={3c^{8/3} (5^{5/3}) (2^{-10/3}}) \eqno(\eqn{CONETWO})
$$
and
$$\eqalign{
	g_1(\tau)&=e^{\beta\tau} \big(R^{1/2}
            e^{\beta\tau/2}+2 K^{1/2}\big)^2, \qquad
	g_2(\tau)=e^{\beta\tau}
            \big(R^{1/2} e^{\beta\tau/2}+2 K^{1/2}\big)^{8/3}.\cr
}
$$

Further define
$$
	M(\tau)=e^{-\nu\lambda\tau}\Big( 1+ \int_0^\tau
		g(s) e^{\nu\lambda s} ds\Big).
	\eqno(\eqn{DEFEMM})
$$
Note that $M(0)=1$.  Differentiating yields
$$
	M'(\tau)=-\nu\lambda M(\tau) + g(\tau).
$$
Therefore
$$
	M'(0)=-\nu\lambda+C_2\Big({\lambda\over \nu^4\lambda_1}\Big)^{1/4}
		\big(R^{1/2}+2K^{1/2}\big)^2
	+C_3\Big({1\over \nu^{5}\lambda_1}\Big)^{1/3}
		\big(R^{1/2}+2K^{1/2}\big)^{8/3}.
$$
Taking $C_4=9c^{8/3}\ge \max\{ (2C_2)^{4/3},2C_3\}=3c^{8/3}(5^{5/3})(2^{-7/3})$ 
we find that $M'(0)<0$
is guaranteed when
$$
	\lambda>{9\over \lambda_1^{1/3}}
		\Big({2c K^{1/2}+c R^{1/2} \over\nu }\Big)^{8/3}.
$$

Now choose $t^*>0$ such that $M(h)<1$ for all $h\in(0,t^*]$.  Note
that $t^*$ only depends on $K$, $\|\delta(t_0)\|$, $\nu$, $\Omega$
and $\lambda$.
We now show the approximating solution $u$ given by (\DASIM) with 
$t_n=hn$ converges to the reference solution $U$ of (\FTWODNS)
as $t\to\infty$.

For induction on $n$ we suppose that $\|\delta(t_n)\|\le R$.
In the case $n=0$ the induction hypothesis is true by definition.
Take inner product of $AQ_\lambda\delta$ with (\DIFEQ) to obtain
$$
	{1\over 2}{d\over dt}\|Q_\lambda\delta\|^2
	+\nu |A Q_\lambda\delta|^2
	\le B_1+B_2+B_3
$$
where
$$B_1=|(B(\delta,U),A Q_\lambda\delta)|,\quad
	B_2=|(B(U,\delta),A Q_\lambda\delta)|\quad\hbox{and}\quad
	B_3=|(B(\delta,\delta),A Q_\lambda\delta)|.
$$

Estimate using \NONLIN\ followed by \PYTHAG\ to obtain
$$\eqalign{
	B_1
	&\le c\lambda_1^{-1/8}\|\delta\|^{3/4}
		|A\delta|^{1/4}\|U\||AQ_\lambda\delta|\cr
	&\le c \Big({2\lambda\over \lambda_1}\Big)^{1/8}
		 \|\delta\|
\|U\||AQ_\lambda\delta|
	+c \Big({2\over\lambda_1}\Big)^{1/8}
		\|\delta\|^{3/4}\|U\| |AQ_\lambda \delta|^{5/4}
\cr
}$$
similarly
$$\eqalign{
	B_2
	&\le c \Big({2\lambda\over \lambda_1}\Big)^{1/8}
		 \|\delta\|
\|U\||AQ_\lambda\delta|
	+c \Big({2\over\lambda_1}\Big)^{1/8}
		\|\delta\|^{3/4}\|U\| |AQ_\lambda \delta|^{5/4}
\cr
}$$
and
$$\eqalign{
	B_3
	&\le c \Big({2\lambda\over \lambda_1}\Big)^{1/8}
		 \|\delta\|^2 |AQ_\lambda\delta|
	+c \Big({2\over\lambda_1}\Big)^{1/8}
		\|\delta\|^{7/4} |AQ_\lambda \delta|^{5/4}.
\cr
}$$
It follows that
$$
	\sum_{i=1}^3 B_i=J_1+J_2
$$
where by \YOUNGS\ with $p=q=2$ we have
$$\eqalign{
	J_1&= c \Big({2\lambda\over \lambda_1}\Big)^{1/8}
	\|\delta\| (\|\delta\|+2\|U\|) |AQ_\lambda\delta|\cr
	&\le {C_2\over 2}\Big({\lambda\over \nu^4\lambda_1}\Big)^{1/4}
		\|\delta\|^2 (\|\delta\|+2\|U\|)^2+
		{\nu\over 4}|AQ_\lambda\delta|^2}
$$
and \YOUNGS\ with $p=8/3$ and $q=8/5$ we have
$$\eqalign{
	J_2&= c \Big({2\over\lambda_1}\Big)^{1/8}
        \|\delta\|^{3/4}(\|\delta\|+2\|U\|) |AQ_\lambda \delta|^{5/4}\cr
	&\le {C_3\over 2} 
		\Big({1\over \nu^5\lambda_1}\Big)^{1/3}
		\|\delta\|^2 (\|\delta\|+2\|U\|)^{8/3}+
		{\nu\over 4}|AQ_\lambda\delta|^2.}
$$
Here, $C_2$ and $C_3$ are as defined in (\CONETWO).
It follows that
$$\eqalign{
	{d\over dt}\|Q_\lambda\delta\|^2
		+\nu|A Q_\lambda\delta|^2
	&\le J_3+J_4}$$
where by \EXGROWTH
$$\eqalign{
		J_3
		&= C_2\Big({\lambda\over \nu^4\lambda_1}\Big)^{1/4}
		\|\delta\|^2 (\|\delta\|+2\|U\|)^2\cr
		&\le C_2\Big({\lambda\over \nu^4\lambda_1}\Big)^{1/4}
			\|\delta(t_n)\|^2 e^{\beta(t-t_n)} \big(\|\delta(t_n)\| 
			e^{\beta(t-t_n)/2}+2\|U\|\big)^2\cr
		&\le C_2\Big({\lambda\over \nu^4\lambda_1}\Big)^{1/4}
			\|\delta(t_n)\|^2 e^{\beta(t-t_n)} \big(R^{1/2}
			e^{\beta(t-t_n)/2}+2 K^{1/2}\big)^2\cr
}
$$
and
$$\eqalign{
     	J_4
		&=C_3 \Big({1\over \nu^5\lambda_1}\Big)^{1/3}
		\|\delta\|^2 (\|\delta\|+2\|U\|)^{8/3}\cr
		&\le C_3 \Big({1\over \nu^5\lambda_1}\Big)^{1/3}
		\|\delta(t_n)\|^2 e^{\beta(t-t_n)} 
			\big(\|\delta(t_n)\| e^{\beta(t-t_n)/2}+2\|U\|\big)^{8/3}\cr
		&\le C_3 \Big({1\over \nu^5\lambda_1}\Big)^{1/3}
		\|\delta(t_n)\|^2 e^{\beta(t-t_n)} 
			\big(R^{1/2} e^{\beta(t-t_n)/2}+2 K^{1/2}\big)^{8/3}.\cr
}
$$

Then by \POINCARE\ and the fact that $P_\lambda\delta(t_n)=0$
we have
$$
	    {d\over dt}\|Q_\lambda\delta\|^2
        +\nu\lambda\|Q_\lambda\delta\|^2
		\le \|Q_\lambda\delta(t_n)\|^2 g(t-t_n)
$$
where $g$ is the function defined in (\DEFGEE).
Multiply by the integrating factor $e^{\nu\lambda(t-t_n)}$
and integrate from $t_n$ to $t$ to obtain
$$
	\|Q_\lambda\delta(t)\|^2
		\le M(t-t_n)\|Q_\lambda\delta(t_n)\|^2
$$
where $M$ is the function defined in (\DEFEMM).  Let
$\gamma=M(h)$.  By our choice of $t^*$ and $h$ we have
$\gamma<1$.  It follows that
$$\eqalign{
	\|\delta(t_{n+1})\|^2&=
	\|Q_\lambda \delta(t_{n+1})\|^2
		\le\lim_{t\nearrow t_{n+1}} \|Q_\lambda \delta(t)\|^2\cr
	&\le\lim_{t\nearrow t_{n+1}} M(t-t_n)\| Q_\lambda\delta(t_n)\|^2
	= M(h)\|  \delta(t_n)\|^2
	= \gamma\| \delta(t_n)\|^2.}
$$
Therefore $\|\delta(t_{n+1})\|^2\le R$, which completes the induction.

To finish the proof note that under these hypothesis we have, 
in fact, proven
$$
	\|\delta(t_n)\|^2\le \gamma^n R.
$$
The proof now finishes as in (\UBOUNDB).
\eprfn

\cor{BIGH} Given any $t^*>0$ there exists $\lambda$ large enough
depending only on
$K$, $\|\delta(t_0)\|$, $\nu$, $\Omega$ and $t^*$ such
that for any $h\in(0,t^*]$ the approximating solution $u$ given
by (\DASIM) with $t_n=hn$ converges to the
reference solution $U$ of (\FTWODNS) as $t\to\infty$.
\prf
First estimate $g(\tau)$ from (\DEFGEE) as
$$\eqalign{
	g(\tau)
	&\le C_2 \Big({\lambda\over \nu^4\lambda_1}\Big)^{1/4}
		L_1^2 e^{2\beta\tau}
	+C_3 \Big({1\over\nu^5\lambda_1}\Big)^{1/3}
		L_1^{8/3} e^{(7\beta/3)\tau}\cr
	&\le \nu \lambda^{1/4} L_2 e^{(7\beta/3)\tau}
}$$
where 
$$L_1=R^{1/2}+2K^{1/2}\quad\hbox{and}\quad
L_2=
    {C_2\over \lambda_1^{1/4}}
        \Big({L_1\over \nu}\Big)^2 
    +{C_3 \over \lambda_1^{7/12}}
        \Big({L_1\over\nu}\Big)^{8/3}.
$$
Therefore
$$\eqalign{
	M(\tau)
	&\le e^{-\nu\lambda\tau}\Big(1+
		\nu\lambda^{1/4} L_2\int_0^\tau e^{(7\beta/3+\nu\lambda)s} ds
		\Big)\cr
	&= e^{-\nu\lambda\tau}\Big(1+
		{\nu\lambda^{1/4} L_2 \over 7\beta/3+\nu\lambda}
			(e^{(7\beta/3+\nu\lambda)\tau}-1)\Big)\le m(\tau)\cr
}$$
where
$$
	m(\tau)= \big(1-L_2\lambda^{-3/4}\big) e^{-\nu\lambda\tau}+
		L_2\lambda^{-3/4} e^{(7\beta/3)\tau}.
$$
Differentiating yields 
$$\eqalign{
	m'(\tau)&=-\nu\lambda \big(1-L_2\lambda^{-3/4}\big) 
			e^{-\nu\lambda\tau}
		+{7\beta\over 3} L_2\lambda^{-3/4}e^{(7\beta/3)\tau}
\quad\hbox{and}\cr
	m''(\tau)&=\nu^2\lambda^2 \big(1-L_2\lambda^{-3/4}\big) 
			e^{-\nu\lambda\tau}
		+{49\beta^2\over 9} L_2\lambda^{-3/4}e^{(7\beta/3)\tau}.
}$$

Given $t^*>0$ choose $\lambda$ large enough such that $m'(t^*)<0$.
Clearly $m''>0$ for this value of $\lambda$.  Thus $m'$
is a strictly increasing function.  It follows that
$m'(s)<0$ for $s\in [0,t^*]$ and therefore
$$
	M(h)\le m(h)=1+\int_0^h m'(s)ds<1$$
for $h\in (0,t^*]$.  Hence, the approximating
solution $u$ given by (\DASIM) with $t_n=hn$ converges to
the reference solution $U$ of (\FTWODNS) as $t\to\infty$.
\eprfn

\cor{ZEROETA} If we take $\eta=0$ and
$$
	\lambda>{9\over \lambda_1^{5/3}}
		\Big({3c |f| \over\nu^2 }\Big)^{8/3},
$$
then the results of \NSDA\ hold where
$t^*$ may be chosen independent of $\|\delta(t_0)\|$.
\prf
When $\eta=0$ then 
$\|\delta(t_0)\|^2=\|QU(t_0)\|^2\le K$.
\eprf

The bound in \ZEROETA\ would be comparable to the bound 
in \NOINDUCT\ if the exponent of $8/3$ were 
instead 2.  The power $8/3$ comes as a result of using the $L^{8/3}$
norm and the particular form of
Agmon's inequality we have used
in estimating the nonlinear term.
Using the same proof technique with different interpolation
inequalities in place of \SOBOLEV, this exponent could be reduced to 
as near $2$ as one might like at the expense of increasing 
the constant $c$.

Next we state without proof the analog of \VARHL\ for the 
incompressible two-dimensional Navier--Stokes equations.

\cor{NSDAEXT}
Let $t^*$ be the bound given in \BIGH.
Suppose $t_{n+1}-t_n\le t^*$ where $t_n\to\infty$ as $n\to\infty$.
Then the approximating solution $u$ given by (\DASIM) converges
to the reference solution $U$ of (\FTWODNS) as $t\to\infty$.
\ecor

We finish with the equivalent of \BNDSLAVE\ for
the Lorenz system.
This result is interesting because it shows that even 
if the approximating solution doesn't converge to the 
reference solution it is still bounded.
This is striking because in the case of 
continuous data assimilation 
it is unknown whether the approximating solution is in
general bounded or not.
In particular, the comments before Theorem 3.5 in 
Olson and Titi [\OLSTIT] indicate that the approximate
solution obtained by continuous data assimilation is not known 
to be bounded if it does not converge to the reference solution.

\thm{UNTA} {\it
There exists $M_5$ independent of $h$ and
depending only on $|f|$, $\Omega$ and $\nu$ such that
$
	\|u(t)\|^2\le M_5 /(1-e^{-\nu\lambda_1 h})$
for all $t$.
}
\prf
Multiply (\FTWODNS) by $Au$, integrate over $\Omega$ and 
apply the inequalities of 
Cauchy--Schwartz and Young to obtain
$$
	{1\over 2} {d\|u\|^2\over dt}
		+\nu |Au|^2 = (f,Au)
		\le |f||Au|
		\le {\nu\over 2}|Au|^2 + {1\over 2\nu} |f|^2.
$$
It follows from \POINCARE\ that 
$$
	{d\|u\|^2\over dt}
		+ \nu\lambda_1 \|u\|^2 \le {1\over \nu} |f|^2.
$$
The rest of the proof is similar to the proof of \BNDSLAVE.
\eprfn

\cor{UBNDNS}
If
$$	\lambda>{9\over \lambda_1^{1/3}}
		\Big({2c K^{1/2}+c\|\delta(t_0)\|\over\nu }\Big)^{8/3}$$
then there exists a bound $M_6$ depending only on
$K$, $\|\delta(t_0)\|$, $\nu$, $\Omega$ and 
$\lambda$ such that the approximating solution $u$ obtained with
$t_n=hn$ is bounded by $M_6$ for any $h>0$.
\prf
The proof is the same as the proof of \UBOUND\ for the Lorenz equations.
Note
that if $\eta=0$ as in \ZEROETA\ then $\lambda$ and $M_6$ may be
chosen independent of $\|\delta(t_0)\|$.
\eprfn

\beginsection{Concluding Remarks}
We have studied discrete data assimilation for the
Lorenz system and the incompressible two-dimensional Navier--Stokes 
equations.  Comparing the results we have obtained for discrete
data assimilation to prior studies of continuous data assimilation 
we find the following.
For the Lorenz system Pecora and Carroll [\PECAR] showed that
continuous data assimilation of the $X$ variable lead to convergence 
of the approximating solution to the reference
solution as time tends to infinity.  In \LORENZDA\ we provide a 
similar result for discrete
data assimilation provided the update time
interval $h$ is sufficiently small.
For the incompressible two-dimensional Navier--Stokes equations 
Olson and Titi [\OLSTIT] obtained conditions on the resolution
parameter $\lambda$ under which the approximating solution
converged to the reference solution.
\NSDA\ states a similar condition on $\lambda$ that leads to
convergence of the approximating solution to the reference solution
provided $h$ is sufficiently small.
In \BIGH\ we also show that for any $h>0$ there is
$\lambda$ large enough such that the approximating 
solution converges to the reference solution as time tends to infinity.
Thus, discrete data assimilation has been shown to work under
similar conditions as continuous data assimilation.

A striking difference between discrete and continuous data assimilation
is given by \BNDSLAVE\ and \UNTA\ which show that the approximating 
solution obtained by discrete data assimilation is bounded even if 
it doesn't converge to the reference solution.  
Although the approximating solution appears bounded in the case of 
continuous data assimilation for all numerical experiments
performed to date, there does not yet exist an analytic proof of
this property.  
Boundedness of an approximating solution that does not converge
to the reference solution remains a conjecture
for continuous data assimilation.

In this work we have found analytic bounds on the update
time interval $h$ for discrete data assimilation which guarantee 
that the approximating solution 
converges to the reference solution.
It is natural to compare our analytic bound with 
numerical simulation.
For the Lorenz system \ANBND\ indicates that
for $t^*\approx 0.000129$ for
the standard parameter values 
$\sigma=10$, $b=8/3$ and $r=28$.
For these same parameter values Hayden~[\KEVIN]
performed a numerical simulation of discrete data assimilation
using the 150-digit-precision variable-step-size 
variable-order Taylor-method integrator [\TAYLOR].
In this work the maximum value for $t^*$ was 
found numerically to lie in the interval $[0.175,0.1875]$.
Thus, convergence of the approximating solution to the reference
solution numerically occurs for values of $h$ three orders
of magnitude larger than those guaranteed by our analysis.

For the incompressible two-dimensional Navier--Stokes equations
bounds on $\lambda$ for continuous data assimilation appear 
in [\OLSTIT] and [\OLSTITB].  That work
shows that
the approximating solution converges to the reference solution
for values of $\lambda$ smaller than expected from the analysis.
As similar techniques are used 
to treat the discrete assimilation
in this paper,  we expect our bounds on $\lambda$ to be 
similarly conservative and that the approximating solution 
obtained by discrete data assimilation will converge numerically 
for much smaller values of $\lambda$ than given by 
\NSDA\ and \ZEROETA.
We also expect our bounds 
on $h$ to be conservative.
A computational study of discrete data assimilation for the 
incompressible two-dimensional Navier--Stokes equations is
currently in progress.

We conclude by returning to the motivating problem of 
using 
satellite imaging data to initialize 
a weather forecasting model.
In applications the value of $h$ 
governing the time interval between consecutive observational 
measurements is generally much larger than the value of 
$\Delta t$ used by the numerical integrator.
Therefore, it is more realistic to treat the observational data 
as measurements occurring at a sequence of times $t_n$ as
was done in this paper
rather than as measurements occurring continuously in time.
In the context of the Lorenz equations and the
incompressible two-dimensional Navier--Stokes equations
we have obtained similar theoretical 
results for discrete data assimilation as for continuous
data assimilation.
We hope that these results will shed light on the differences and
similarities between discrete and continuous data assimilation and
help guide future work in understanding more complicated problems.

\beginsection{Acknowledgements}
The work of E.S.T. was supported in part by NSF grants
no.~DMS--0708832 and DMS-1009950, and by
the Alexander von Humboldt Stiftung/Foundation and the 
Minerva Stiftung/Foundation.

\vfill\eject
\beginsection{References}
{\frenchspacing\parskip=5pt plus2pt minus1pt

\ref{BREZIS}
H. Br\'ezis, T. Gallouet,
Nonlinear Schr\"odinger evolution equations,
{\it Nonlinear Anal.} {\bf 4} (1980), no. 4, 677--681.

\ref{CONFOI}
P. Constantin, C. Foias,
{\it Navier--Stokes equations},
Chicago Lectures in Mathematics, 1988.

\ref{CONFOITEM}
P. Constantin, C. Foias, R. Temam,
Attractors Representing Turbulent Flows,
{\it Mem. Amer. Math. Soc.} {\bf 53} (1985), no. 314.

\ref{FJKT}
C. Foias, M.S. Jolly, I. Kukavica, E.S Titi,
The Lorenz Equation as a Metaphor for the Navier-Stokes Equations,
{\it Discrete and Continuous Dynamical Systems},
{\bf 7} (2001), no. 4, 403--429.

\ref{FOIPRO}
C. Foias, G. Prodi,
Sur le comportement global des solutions non-stationnaires des
e\'quations de Navier--Stokes en dimension 2,
{\it Rend. Sem. Mat. Univ. Padova} {\bf 39} (1967), 1-34.

\ref{KEVIN}
K. Hayden,
{\it Synchronization in the Lorenz System}, 
Masters Thesis, University of Nevada, 
Department of Mathematics and Statistics, 2007.

\ref{JONTITA}
D.A.~Jones, E.S.~Titi,
On the Number of Determining Nodes for the 2D
Navier--Stokes Equations,
{\it Journal of Mathematical Analysis and Applications},
{\bf 168} (1992), no. 1, 72--88.

\ref{JONTITB}
D.A.~Jones, E.S.~Titi,
Upper Bounds on the Number of Determining Modes,
Nodes, and Volume Elements for the Navier--Stokes Equations,
{\it Indiana Univ. Math. J.}, {\bf 42} (1993), no. 3, 875--887.

\ref{TAYLOR}
A.~Jorba, M.~Zou,
A Software Package for the Numerical Integration of ODEs
by Means of High-Order Taylor Methods,
{\it Experimental Mathematics}, {\bf 14} (2005), no. 1, 99--117.

\ref{LADY}
O.A.~Ladyzhenskaya,
{\it The Mathematical Theory of Viscous
Incompressible Flow},
Translated from the Russian by R.~Silverman
and John Chu,
Gordon and Breach Science Publishers, 1969.

\ref{OLSTIT}
E. Olson, E.S. Titi,
Determining Modes for Continuous Data Assimilation in 2D
Turbulence, {\it Journal of Statistical Physics} {\bf 113} (2003),
no. 516, 799--840.

\ref{OLSTITB}
E. Olson, E.S. Titi,
Determining Modes and Grashoff Number for Continuous Data Assimilation in 2D
Turbulence, {\it Theoretical and Computational Fluid Dynamics}, {\bf 22}
(2008), 327-339.

\ref{PSEUDO}
E. Olson, E.S. Titi,
Viscosity Versus Vorticity Stretching: Global Well-posedness
for a family of Navier--Stokes Alpha-like Models,
{\it Nonlinear Analysis}, {\bf 66} (2007), no. 5, 1635--1673.

\ref{PECAR}
L. Pecora, T. Carroll,
Synchronization in Chaotic Systems,
{\it Physical Review Letters}, {\bf 64} (1990), no. 8, 821--824.

\ref{ROB}
J.C. Robinson,
{\it Infinite-Dimensional Dynamical Systems:
An Introduction to Dissipative Parabolic PDEs and the
Theory of Global Attractors},
Cambridge Texts in Applied Mathematics, 2001.

\ref{TEMA}
R. Temam,
{\it Infinite-dimensional Dynamical Systems in Mechanics
and Physics, 2nd Edition}, Springer--Verlag, New York, 1997.

\ref{TEMB}
R. Temam,
{\it Navier--Stokes Equations and Nonlinear Functional Analysis},
second edition,
CBMS-NSF Regional Conference Series in Applied Mathematics vol 66,
SIAM, 1995.

\ref{WINGARD}
C. Wingard, {\it Removing Bias and Periodic Noise in
Measurements of the Lorenz System},
Thesis, University of Nevada, 
Department of Mathematics and Statistics, 2009.

\ref{WU} X. Wu, J. Lu, C. Wang, J. Liu,
Impulsive Control and Synchronization of the Lorenz Systems Family,
{\it Chaos, Solitons and Fractals}, {\bf 31} (2007), no. 3, 631--638.

\ref{YANG} T. Yang, L. Yang, C. Yang,
Impulsive Control of Lorenz System,
{\it Physica D}, {\bf 110} (1997), no. 1-2, 18--24.

}

\vfill\end